\documentclass[12pt,reqno,a4wide]{amsart}

\usepackage{tikz} 

\usetikzlibrary{decorations.pathreplacing}
\usetikzlibrary{arrows,decorations.pathmorphing,snakes}
\usetikzlibrary{fit}
\usepackage{calrsfs}
\usepackage[charter]{mathdesign}

\allowdisplaybreaks

\title[Counting ternary trees according to the number of middle edges]{Counting ternary trees according to the number of middle edges and factorizing into $(3/2)$-ary trees}

\author{Helmut Prodinger}
\address{Helmut Prodinger\\
	Department of Mathematical Sciences\\
	Stellenbosch University\\
	7602 Stellenbosch,	South Africa, and
	 Department of Mathematics and Mathematical Statistics\\
	 Umea University\\
	 907 36 Umea, 	 Sweden  }
\email{hproding@sun.ac.za}

\keywords{Ternary trees, middle edges, cubic equation}
\subjclass[2020]{05A15, 05A16}

\begin{document}
	
	\begin{abstract}
		The sequence A120986 in the Encyclopedia of Integer Sequences counts ternary trees according to the number of nodes and the number of middle edges. Using a certain substition, the underlying cubic equation can be factored. This leads to an extension
		of the concept of $(3/2)$-ary trees, introduced by Knuth in his christmas lecture from 2014.
	\end{abstract}
	
	\maketitle

	\section{Introduction}
	The recent preprint \cite{Burstein} triggered my interest in the sequence A120986 in \cite{OEIS}. The double-indexed sequence
	enumerates ternary trees according to the number of edges and the number of middle edges. We consider here $T(n,k)$, the
	number of ternary trees with $n$ nodes and $k$ middle edges. The difference is marginal, but we want to compare/relate our analysis with \cite{christmas}, and there it is also the number of nodes that is considered. Let
	$G=G(x,u)=\sum_{n,k\ge0}T(n,k)x^nu^k$. Then it is easy to see (decomposition at the rooot) that
	\begin{equation*}
G=1+xG^2(1-u+uG).
	\end{equation*}
The substitution
\begin{equation*}
x=\frac{t(1-t)^2}{(1-t+ut)}
\end{equation*}
makes the cubic equation manageable and also allows, as in \cite{christmas}, to introduce a (refined) version of the $(3/2)$-ary trees.

Here is a small table of these numbers and a ternary tree:
\begin{center}
	
	\begin{tabular}{c|ccccccccccc}
		
		$n\backslash k$& 0&1&2&3&4&5\\
		\hline
	0&1&&&&&&&\\
		1&1&&&&&&&\\
		2&2&1&&&&&&\\
		3&5&6&1&&&&&\\
		4&14&28&12&1&&&&\\
		5&42&120&90&20&1&&&\\
		6&132&495&550&220&30&1&&\\
	\end{tabular}
	
\end{center}
\begin{figure}[h]
	\begin{tikzpicture}[scale=0.5]
	
	\draw[](0,0)to(4,-4);
	
	\node at (0,0) {$\bullet$};\node at (1,-1) {$\bullet$};\node at (2,-2) {$\bullet$};\node at (3,-3) {$\bullet$};
	\node at (4,-4) {$\bullet$};
	
	\draw(0,0)to(-3,-3);\node at (-1,-1) {$\bullet$};\node at (-2,-2) {$\bullet$};\node at (-3,-3) {$\bullet$};

	\draw[ultra thick](3,-3)to(3,-4);
	\draw[ultra thick](1,-1)to(1,-2);\draw(1,-1)to(-3,-5);
	\draw[ultra thick](0,0)to(0,-1);\node at (0,-1) {$\bullet$};
	
	\node at (3,-4) {$\bullet$};
	\node at (1,-2) {$\bullet$};\node at (0,-2) {$\bullet$};\node at (-1,-3) {$\bullet$};
	\node at (-2,-4) {$\bullet$};\node at (-3,-5) {$\bullet$};
	
	\draw(3,-3)to(1,-5);\node at (1,-5) {$\bullet$};\node at (2,-4) {$\bullet$};

	\end{tikzpicture}
	
	\caption{Ternary tree with 17 nodes and 3 middle edges}
	
\end{figure}
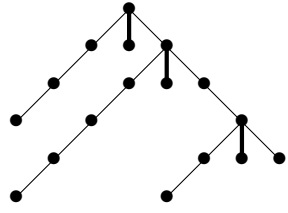

\section{Analysis of the cubic equation}

The cubic equation has the following solutions:
\begin{align*}
r_1&=\frac{1}{1-t},\\
r_2&=\frac{-t+t^2-t^2u+\sqrt{t(1-t+ut)(4u+t-4ut-t^2+t^2u)}}{2ut(1-t)},\\
r_3&=\frac{-t+t^2-t^2u-\sqrt{t(1-t+ut)(4u+t-4ut-t^2+t^2u)}}{2ut(1-t)}.
\end{align*}
Note that
\begin{equation*}
r_2r_3=-\frac{1-t+ut}{ut(1-t)}.
\end{equation*}
The root with the combinatorial significance is $r_1$. But it is the explicit form of the two other roots that makes everything here interesting and challenging. 

We extract coefficients of $r_1$ using contour integration, which is closely related to the Lagrange inversion formula. The path of integration is a small circle in the $x$-plane which is then transformed into a small circle in the $t$-plane.
\begin{align*}
[x^n]r_1&=\frac1{2\pi i}\oint \frac{dx}{x^{n+1}}\frac1{1-t}\\
&=\frac1{2\pi i}\oint \frac{dt(1-t)(1-3t+2t^2-2t^2u)}{(1-t+tu)^2}\frac{(1-t+tu)^{n+1}} {t^{n+1}(1-t)^{2n+2}}\frac1{1-t}\\
&=[t^n](1-3t+2t^2-2t^2u)\frac{(1-t+tu)^{n-1}} {(1-t)^{2n+2}}.
\end{align*}
Furthermore
\begin{align*}
	[x^nu^k]r_1
	&=[t^n][u^k](1-3t+2t^2-2t^2u)\frac{(1-t+tu)^{n-1}} {(1-t)^{2n+2}}\\
	&=[t^n]\binom{n-1}{k}\frac{t^k(1-2t)}{(1-t)^{n+k+2}}	-2[t^n]\binom{n-1}{k-1}\frac{t^{k+1}}{(1-t)^{n+k+2}}\\
		&=\binom{n-1}{k}[t^{n-k}]\frac{(1-2t)}{(1-t)^{n+k+2}}	-2\binom{n-1}{k-1}[t^{n-k-1}]\frac{1}{(1-t)^{n+k+2}}\\
		&=\binom{n-1}{k}\binom{2n+1}{n-k}-2\binom{n-1}{k}\binom{2n}{n-k-1}
		-2\binom{n-1}{k-1}\binom{2n}{n-k-1}\\
		&=\frac1n\binom nk\binom{2n}{n-1-k}.
\end{align*}
For $u=1$, which means that the middle edges are not especially counted, we get
\begin{equation*}
\sum_k\frac1n\binom nk\binom{2n}{n-1-k}=\frac1n\binom{3n}{n-1},
\end{equation*}
the number of ternary trees with $n$ nodes.

\section{Factorizing the solution of the cubic equation}

For $u=1$, Knuth \cite{christmas} was able to factor the generating function $r_1$ into two factors, for which he coined the catchy name $(3/2)$-ary trees. For this factorization, see also \cite{naimi-paper,BM-P}. The goal in this section is to perform this factorization in the context of counting middle edges, i. e., for the generating function with the additional variable $u$. In Knuth's instance, the generating function was expressible as a generalized binomial series (in the sense of Lambert \cite{GKP}), but that does not seem to be an option here.

Note that
\begin{equation*}
\frac1{r_2}=\frac t2-\frac{\sqrt t\sqrt{t(1-t)+u(2-t)^2}}{2\sqrt{1-t+tu}}
\end{equation*}
and
\begin{equation*}
	\frac1{r_3}=\frac t2+\frac{\sqrt t\sqrt{t(1-t)+u(2-t)^2}}{2\sqrt{1-t+tu}}.
\end{equation*}
From the cubic equation we deduce that
\begin{equation*}
r_1=-\frac1{uxr_2r_3},
\end{equation*}
which is the desired factorization. The factor $ux$ will be fairly split as $\sqrt{ux}\cdot\sqrt{ux}$, whereas the minus sign goes to the
factor $1/r_2$. In the following we work out how this factorization can be obtained. To say it again, it is not as appealing as in the original case.

Let us write
\begin{equation*}
t=x\Phi(t), \quad\text{with}\quad  \Phi(t)=\frac{1-t+tu}{(1-t)^2},
\end{equation*}
so that we can use the Lagrange inversion formula to get
\begin{equation*}
[x^n]t^{\ell}=\frac{\ell}{n}[t^{n-\ell}]\Phi(t)^n
\end{equation*}
and
\begin{align*}
	[x^nu^k]t^{\ell}&=\frac{\ell}{n}[t^{n-\ell}][u^k]\frac{(1-t+tu)^n}{(1-t)^{2n}}\\
	&=\frac{\ell}{n}[t^{n-\ell-k}]\binom{n}{k}\frac1{(1-t)^{n+k}}
	=\frac{\ell}{n}\binom{n}{k}\binom{2n-\ell-1}{n-\ell-k}.
\end{align*}
In particular,
\begin{equation*}
t=\sum_{n\ge1}x^n\sum_{0\le k\le n}\frac{1}{n}\binom{n}{k}\binom{2n-2}{n-1-k}u^k;
\end{equation*}
this series expansion may be used in the following developments whenever needed.

To proceed further, we set $u=1+U$ and $\tau=t/u$:
\begin{gather*}
\frac1{r_2}=\frac t2-\frac{\sqrt{x}}{2(1-t)}\sqrt{4-3t+U(2-t)^2},\\
	\frac1{r_3}=\frac t2+\frac{\sqrt{x}}{2(1-t)}\sqrt{4-3t+U(2-t)^2}
\end{gather*}
Since the first term is well understood, we concentrate on the second:
\begin{multline*}
\frac{\sqrt{x}}{2(1-t)}\sqrt{4-3t+U(2-t)^2}
\\=\sqrt{ux}\bigg(	1+\frac18 \left( 5+4U \right) \tau+\frac {1}{128}  
( 71+136U	+64{U}^{2} ) {\tau}^{2}\\+\frac {1}{1024} ( 541+1596U	+1568{U}^{2}+512{U}^{3} ) {\tau}^{3}+\cdots\bigg)	=:\sqrt{ux}\cdot\Xi.
	\end{multline*} 
	With this expanded form $\Xi$, we have now our final formula, the expansion of $r_1$ into two factors:
\begin{align*}
r_1=\frac1{1-t}=-\frac{1}{ux}\frac1{r_2}\frac1{r_3}
=\Big(-\frac{1}{\sqrt{ux}}\frac t2+\Xi\Big)\Big(\frac{1}{\sqrt{ux}}\frac t2+\Xi\Big).
\end{align*}
These two factors do not have a combinatorial meaning, as it seems, but we can still stick to the $(3/2)$-ary tree notation, with the additional counting of middle edges.

\end{document}